\documentclass[draft]{amsart}
\title{
  The duality of conformally flat manifolds
}

\date{June 14, 2010}
\usepackage[usenames]{color}
\author{Huili~Liu}
\address[Liu]{%
 Department of Mathematics,
  Northeastern University,
  Shenyang 110004,  PR China}
\email{liuhl@mail.neu.edu.cn}

\author{Masaaki Umehara}
\address[Umehara]{%
   Department of Mathematics, Graduate School of Science,
   Osaka University,
   Toyonaka, Osaka 560-0043,
   Japan
}
\email{umehara@math.sci.osaka-u.ac.jp}

\author{Kotaro Yamada}
\address[Yamada]{%
   Department of Mathematics,
   Tokyo Institute of Technology,
   O-okayama, Meguro, Tokyo 152-8551, Japan}
\email{kotaro@math.titech.ac.jp}

\thanks{%
 The second and the third authors were partially 
 supported by Grant-in-Aid for Scientific Research 
 (A) No. 19204005 and (B) No. 21340016,
 respectively from the Japan Society for the Promotion of Science.
}
\newtheorem{thm}{Theorem}[section]
\newtheorem{prop}[thm]{Proposition}
\newtheorem{lem}[thm]{Lemma}
\newtheorem{fact}[thm]{Fact}
\newtheorem{cor}[thm]{Corollary}
\theoremstyle{definition}
\newtheorem*{defi}{Definition}
\newtheorem{ex}[thm]{Example}
\newtheorem{rmk}[thm]{Remark}

\newtheorem*{acknowledgement}{Acknowledgement}

\numberwithin{equation}{section}

\newcommand{\Ric}{\operatorname{Ric}}
\newcommand{\Cod}{\operatorname{Cod}}
\newcommand{\Hess}{\operatorname{Hess}}
\newcommand{\trace}{\operatorname{Trace}}

\newcommand{\R}{\boldsymbol{R}}
\newcommand{\Lor}{\boldsymbol{L}}
\newcommand{\C}{\boldsymbol{C}}
\newcommand{\GCF}{\mathcal G \hskip -1pt C\mathcal\!\mathcal F}
\newcommand{\MFr}{\mathcal M_{Fr}}
\newcommand{\Reg}{\mathcal R}

\renewcommand{\O}{\mathcal O}
\newcommand{\E}{\mathcal E}
\newcommand{\D}{\mathcal D}
\newcommand{\lto}{\longrightarrow}
\newcommand{\inner}[2]{\left\langle{#1},{#2}\right\rangle}

\renewcommand{\phi}{\varphi}
\renewcommand{\epsilon}{\varepsilon}
\newcommand{\pmt}[1]{{\begin{pmatrix} #1  \end{pmatrix}}}

\newcommand{\first}{\operatorname{\mathit{I}}}
\newcommand{\second}{\operatorname{\mathit{I\!I}}}
\newcommand{\third}{\operatorname{\mathit{I\!I\!I}}}
\newcommand{\Ker}{\operatorname{Ker}}

\begin{document}
\maketitle

\begin{abstract}
In a joint work with Saji, the second and the third authors gave
an intrinsic formulation of wave fronts and proved a 
realization theorem for wave fronts in space forms.
As an application, we show that the following four objects
are essentially the same:
 \begin{itemize}
  \item conformally flat $n$-manifolds $(n\ge 3)$ 
	with admissible singular points
	(i.e. admissible GCF-manifolds),
  \item frontals as hypersurfaces in the lightcone $Q^{n+1}_+$,
  \item frontals as hypersurfaces in the hyperbolic space
	$H^{n+1}$,
  \item spacelike frontals as hypersurfaces in the 
	de Sitter space $S^{n+1}_1$.
\end{itemize}
 Recently, the duality of conformally flat Riemannian manifolds
 was discovered by  several geometers.
 In our setting, this duality can be explained via the existence 
 of a two-fold map of the congruence classes of admissible GCF-manifolds
 into that of frontals in  $H^{n+1}$.
 It should be remarked that the dual conformally flat metric
 may have degenerate points even when the original conformally flat
 metric is positive definite.
 This is the reason why we consider conformally flat manifolds with
 singular points.
 In fact, the duality is an involution on the set of admissible
 GCF-manifolds.
 The case $n=2$ requires a special treatment, since
 any Riemannian $2$-manifold is conformally flat. 
 At the end of this paper, we also determine the moduli space 
 of isometric immersions of  a given simply connected 
 Riemannian $2$-manifold into 
the lightcone $Q^{3}_+$.
\end{abstract}

%%%%%%%%%%%%%%%%%%%%%%%%%%%%%%%%%%%%%%%%%%%%%%%%
\section{Introduction}
We denote by $Q^{n+1}_+$ the upper lightcone in 
Lorentz-Minkowski space $\Lor^{n+2}$.
Izumiya \cite{I} pointed out that the mean curvature $H$ of a surface in
$Q^{3}_+$ is equal to $-1/2$ times of the Gaussian curvature $K$,
as Theorema Egregium of Gauss for surface theory in  $Q^{3}_+$.
(It should be remarked that our notation is different from
that in \cite{I}.
If one use Izumiya's notation, $H$ coincides with $K$.
This formula was found independently by the first author
in \cite[(2.7)]{L} 
and in more general form in \cite[(2.13)]{LJ}.)
When $n\ge 3$, this corresponds to the fact that a conformally flat
$n$-manifold can be isometrically immersed as a hypersurface in 
$Q^{n+1}_+$, and the second fundamental form is just
its Schouten tensor (see Brinkmann \cite{Br} and
Asperti-Dajczer \cite{AD}).

Recently, Izumiya \cite{I}, Espinar-G\'alvez-Mira \cite{EGM} and
Liu-Jung \cite{LJ} independently found a duality of
hypersurfaces in $Q^{n+1}_+$.
More precisely, Izumiya \cite{I} explained this duality on hypersurfaces
in $Q^{n+1}_+$ as a bi-Legendrian fibration in contact geometry.
The two distinct explicit formulas of the dual in the 
lightcone are given in \cite{LJ} and \cite[p.~332]{I} respectively.
On the other hand, Espinar, G\'alvez and Mira \cite{EGM}
found the duality on conformally flat manifolds
from the viewpoint of hypersurface theory in the hyperbolic space,
and found that the inverses of the eigenvalues of their Schouten
tensors coincide with the eigenvalues of the dual Schouten tensors.
It should be remarked that the dual metric of
a conformally flat Riemannian metric might degenerate,
in general.

On the other hand, 
in a joint work \cite{SUY5} with Saji, the second and the third
authors gave the definition of a frontal bundle, and proved 
a realization of it as a wave front in space forms, which is a
generalization of the fundamental theorem of surface theory.
In this paper, as an application of this, we define 
`admissible generalized conformally flat manifolds'
(or `admissible GCF-manifolds') 
as a class of conformally flat manifolds with admissible singular points,
and show that the above duality operation is 
an involution on this class.
Also, we  give an explicit formula for dual metrics, and remark that
the Schouten tensors are invariant under the duality operation.
Moreover, as a refinement of the 
result in \cite{EGM}, under the assumption that $M^n$ $(n\ge 3)$
is $1$-connected (i.e. connected and simply connected), 
we show that this duality 
comes from the existence of the two-fold map 
\[
    \Psi:\GCF(M^n)\longrightarrow \MFr(M^n,H^{n+1})
\]
of the moduli space $\GCF(M^n)$ of admissible GCF-manifolds 
into the moduli space 
$\MFr(M^n,H^{n+1})$ of frontals in hyperbolic space $H^{n+1}$.
To prove the existence of the map $\Psi$,
we apply the realization theorem of intrinsic wave fronts given in
\cite{SUY5}.  

Finally, we consider the $2$-dimensional case, and determine the moduli
space of isometric immersions of a given simply connected 
Riemannian $2$-manifold into $Q^3_+$.

\section{The duality of conformally flat manifolds}
\label{sec:duality}
A Riemannian $n$-manifold $(M^n,g)$ is called 
{\em conformally flat\/}
if for each point $p\in M^n$, there exists a neighborhood
$U(\subset M^n)$ of $p$ and a $C^\infty$-function $\sigma$
on $U$ such that $e^{2\sigma}g$ is a metric with vanishing
sectional curvature.
When $n\ge 4$, $(M^n,g)$ is conformally flat
if and only if
the conformal curvature tensor 
\begin{equation}\label{eq:weyl}
 W_{ijkl}:=R_{ijkl}+
   A_{ik}g_{jl}-A_{il}g_{jk}+A_{jl}g_{ik}-A_{jk}g_{il}
\end{equation}
vanishes identically on $M^n$, where
$(u^1,\dots,u^n)$ is a local coordinate system of $M^n$, 
\begin{equation}\label{eq:schouten}
   A:=\frac{1}{n-2}\sum_{i,j}\left(R_{ij}-\frac{S_g}{2(n-1)}g_{ij}\right)\,
   du^i\otimes du^j
\end{equation}
is called the {\em Schouten tensor}, 
$g_{ij}$, $R_{ijkl}$, $R_{ij}$ are the components of 
the metric $g$, the curvature tensor of $g$,
and the Ricci tensor of $g$ respectively, 
and $S_g$ denotes the scalar curvature.
When $n=3$, $(M^3,g)$ is  conformally flat
if and only if 
$A$ in \eqref{eq:schouten}  is a Codazzi tensor, that is, 
$\nabla A$ is a symmetric $3$-tensor, where $\nabla$ is the Levi-Civita
connection of $(M^3,g)$.
(When $n\ge 4$, conformal flatness implies that
$A$ is a Codazzi tensor because of the second Bianchi identity.)
When $n=2$, all Riemannian metrics are conformally flat.

To formulate conformally flat manifolds with singularities, 
we need to define the following:
\begin{defi}[\cite{SUY1}, \cite{SUY2}, \cite{SUY5}]
 Let $\E$ be a vector bundle over an $n$-manifold $M^n$
 $(n\ge 1)$ of rank $n$, 
 and $\phi:TM^n\to \E$ a bundle homomorphism,
 where $TM^n$ is the tangent bundle of $M^n$.
 Suppose that $\E$ has a metric $\inner{~}{~}$ and
 a metric connection $D$.
 Then $(M^n, \E,\inner{~}{~},D,\phi)$ is called a 
 {\em coherent tangent bundle\/} if $\phi$ satisfies the condition 
 \begin{equation}\label{eq:c}
    D_{X}\phi(Y)-D_{Y}\phi(X)-\phi([X,Y])=0 
 \end{equation}
 for any $C^\infty$-vector fields $X,Y$ on $M^n$.
 Let $(M^n, \E,\inner{~}{~},D,\phi)$ be a coherent tangent
 bundle. 
 Then $p\in M^n$ is called a {\em singular point\/} of
 $\phi$ if the linear map $\phi_p:T_pM^n\to \E_p$ is not injective,
 where $\E_p$ is the fiber of $\E$ at $p$.
 On the other hand, $p\in M^n$ is called a {\it regular point\/}
 if it is not a singular point. We denote by
$\Reg_{M^n}$ or $\Reg_{M^n,\phi}$  the set of regular 
points of $\phi$.
\end{defi}
Coherent tangent bundles can be considered as a generalization of
Riemannian metrics.
In fact, the pull-back metric 
$g:=\phi^*\inner{~}{~}$ gives a 
Riemannian metric on $\Reg_{M^n}$ and the
Levi-Civita connection $\nabla$ of $g$
coincides with the pull-back of the connection $D$ by $\phi$
because of the condition \eqref{eq:c}.
Moreover, one can prove the Gauss-Bonnet formula
when $n=2$ (see \cite{SUY1}, \cite{SUY2} and \cite{SUY5}).
\begin{ex}\label{ex:front}
 Let $M^n$ and $N^{n+1}$ be $C^\infty$-manifolds of dimension 
 $n$ and of dimension $n+1$, respectively.
 The projectified cotangent bundle $P(T^*N^{n+1})$
 has a canonical contact structure. 
 A $C^\infty$-map $f:M^n \to N^{n+1}$ is called a 
 {\em frontal\/} if $f$ lifts to a Legendrian map $L_f$, 
 i.e., a $C^\infty$-map $L_f:M^n\to P(T^*N^{n+1})$
 such that the image $dL_f(TM^n)$ of the tangent
 bundle $TM^n$ lies in the contact hyperplane field on 
 $P(T^*N^{n+1})$. 
 Moreover, $f$ is called
 a {\em wave front\/} or a {\it front\/} if it lifts to a 
 Legendrian immersion $L_f$. 
 Frontals (and therefore fronts) generalize immersions, 
 as they allow for singular points. 
 A frontal  $f$ is said to be co-orientable if its Legendrian 
 lift $L_f$ can lift up to a $C^\infty$-map into
 the cotangent bundle $T^*N^{n+1}$.

 Now, we fix a Riemannian metric 
 $\tilde g$  on $N^{n+1}$.
 Then, it can be easily checked that a $C^\infty$-map
 \[
    f:M^n\longrightarrow N^{n+1}
 \]
 is a frontal if and only if for each $p\in M^n$,
 there exists a neighborhood $U$ of $p$ and a unit 
 $C^\infty$-vector field $\nu$ of $N^{n+1}$ along 
 $f$ defined on $U$ such that $\tilde g\bigl(df(X),\nu\bigr)=0$ 
 holds for any vector fields $X$ on $U$
 (that is, $\nu$ is a locally defined unit normal vector field).
 Moreover, if the locally defined unit normal vector field
 $\nu:U\to T_1N^{n+1}$ can be taken to be an immersion 
 for each $p\in M^n$,  $f$ is called a {\em front}, 
 where $T_1N^{n+1}$ is the unit tangent bundle 
 of $(N^{n+1},\tilde g)$.
 We denote by $\E_f$ the subbundle of the pull-back
 bundle $f^*\bigl(TN^{n+1}\bigr)$
 consisting of vectors perpendicular to $\nu$. 
 Then
 \[
    \phi_f\colon{}TM^n\ni X \longmapsto df(X)\in \E_f
 \]
 gives a bundle homomorphism.
 Let $\tilde\nabla$ be the Levi-Civita connection on $N^{n+1}$.
 Then by taking the tangential part of $\tilde\nabla$, 
 a connection $D$ on $\E_f$ satisfying  \eqref{eq:c} is induced.
 Let $\inner{~}{~}$ be a metric on $\E_f$ induced from 
 the Riemannian metric on $N^{n+1}$, then  $D$ is a metric connection
 on $\E_f$. 
 Thus we get a coherent tangent bundle 
 $(M^n,\E_f,\inner{~}{~},D,\phi_f)$.
\end{ex}
In this setting, we define the following
\begin{defi}
 A given coherent tangent bundle $(M^n, \E,\inner{~}{~},D,\phi)$ 
 $(n\ge 3)$ is called a {\em generalized conformally flat manifold\/}
 or a {\em GCF-manifold\/}  if the regular set 
 $\Reg_{M^n}$ of $\phi$ is dense in $M^n$ and
 the pull-back metric $g:=\phi^*\inner{~}{~}$ is conformally flat
 on $\Reg_{M^n}$.
\end{defi}
\begin{ex}\label{ex:map}
 Let $(N^n,\tilde g)$ $(n\ge 3)$ be a conformally 
flat Riemannian $n$-manifold,
 and $f:M^n\to N^n$ a $C^\infty$-map
 whose regular set is dense in $M^n$. 
 Let $\E_f:=f^*\bigl(TN^n\bigr)$ be the pull-back bundle of 
 $TN^n$ by $f$.
 Then $\tilde g$ induces a positive definite metric 
 $\inner{~}{~}$ on $\E_f$,
 and the pull-back $D$ of the Levi-Civita connection of $\tilde g$
 gives a metric connection on $(\E_f,\inner{~}{~})$.
 We set 
 \[
   \phi^{}_f:=df:TM^n\longrightarrow \E_f.
 \]
 Then we have a coherent tangent bundle 
 $(M^n,\E_f,\inner{~}{~},D,\phi_f)$
 which is a GCF-manifold because $(N^n,\tilde  g)$ is conformally flat.
\end{ex}

\begin{defi}
 Let $(M^n,\E,\inner{~}{~},D,\phi)$ $(n\ge 3)$ be a GCF-manifold.
 A point $p\in \Reg_{M^n}$ is called a {\it parabolic point\/}
 if the Schouten tensor $A$ of the induced metric $g=\phi^*\inner{~}{~}$
degenerates. 
 We denote by $\Reg_{M^n}^*(\subset \Reg_{M^n})$ 
the set of non-parabolic points.
\end{defi}

The following assertion is the explicit description of
the duality of conformally flat manifolds:
\begin{thm}\label{thm:main1}
 Let $(M^n,\E,\inner{~}{~},D,\phi)$ $(n\ge 3)$ be a GCF-manifold,
 and $g=\phi^*\inner{~}{~}$ the induced metric.
 Then 
 \begin{equation}\label{eq:tilde}
  \check g:=\sum_{i,j,a,b} A_{ia}A_{jb}g^{ab}\,du^i\otimes du^j
 \end{equation}
 gives a conformally flat metric on $\Reg^*_{M^n}$
 {\rm(}called the {\em dual metric}{\rm)}, where
 $(u^1,\dots,u^n)$ is a local coordinate system,
 $(g^{ab})$ is the inverse matrix of $(g_{ij})$  and 
 $A$ is  the Schouten tensor of $g$.
 Moreover,  the Schouten tensor of $\check g$ 
 coincides with $A$.
\end{thm}
One can prove the assertion by a direct calculation.
We give an alternative proof in Section~\ref{sec:hyperbolic}.
The following assertion follows immediately, which was
proved in \cite{EGM}:
\begin{cor}[\cite{EGM}]
 The eigenvalues of the Schouten tensor with respect to $\check g$
 are the inverses of those of $A$ with respect to $g$.
\end{cor}

As seen in Example \ref{ex:map}, the class of GCF-manifolds
might be too wide.  
We now define a subclass of GCF-manifolds which admits 
the duality and also the conformal 
changes of metrics as follows:
Set
 \begin{equation}\label{hatA}
  \hat A:=\sum_{i,j,a} g^{ia}A_{aj}\,\frac{\partial}{\partial u^i}\otimes du^j,
 \end{equation}
which is a tensor defined on $\Reg_{M^n}$,
 where $g=\phi^*\inner{~}{~}$ is the induced metric  and 
 $A$ is  the Schouten tensor of $g$ on $\Reg_{M^n}$.

\begin{defi}
 Let $(M^n,\E,\inner{~}{~},D,\phi)$ $(n\ge 3)$ be a
 GCF-manifold. 
 Then it is called {\em admissible\/} if $\Reg_{M^n}^*$ is 
 dense in  $M^n$ and the tensor $\hat A$
 induces a new bundle homomorphism
 \begin{equation}\label{eq:tilde-phi}
   \check \phi:TM^n\ni v \longmapsto \phi\circ \hat A(v)\in \E,
 \end{equation}
 namely, the homomorphism 
 $TM^n|_{\Reg_{M^n}}\ni v\mapsto \phi\circ \hat A(v) \in \E$ 
 can be smoothly extended to the whole of $TM^n$. 
\end{defi}
In this setting, we can formulate the duality of conformally flat
manifolds as follows:
\begin{thm}\label{thm:main2}
 Let $(M^n,\E,\inner{~}{~},D,\phi)$ $(n\ge 3)$  be an admissible
 GCF-manifold.
 Then by replacing $\phi$ by $\check \phi$ in \eqref{eq:tilde-phi},
 $(M^n,\E,\inner{~}{~},D,\check\phi)$
 is also an admissible GCF-manifold.
 Moreover $\check{\check \phi}$ coincides with $\phi$.
\end{thm}
We prove this assertion in Section~\ref{sec:hyperbolic}.
\begin{rmk}
 We can prove Theorems \ref{thm:main1} and \ref{thm:main2}
 under the assumption that $\inner{~}{~}$ is a  
 non-degenerate symmetric bilinear form. 
 So, for example, the duality also holds
 for Lorentzian conformally flat manifolds. 
\end{rmk}
Moreover, we also prove the following in 
Section \ref{sec:frontal}:
\begin{prop} \label{prop:conformal}
 Let $(M^n,\E,\inner{~}{~},D,\phi)$ $(n\ge 3)$ be 
 an admissible GCF-manifold. 
 Then for a given  $C^\infty$-function 
 $\sigma$ on $M^n$, there exist a connection $D^\sigma$ 
 and  a bundle homomorphism
 $\phi_\sigma:TM^n\to \E$ 
 such that $(M^n,\E,e^{2\sigma}\inner{~}{~},D^{\sigma},\phi_{\sigma})$
 is a GCF-manifold.
\end{prop} 

The conformal change of the metric $g$ of 
a GCF-manifold as in Proposition \ref{prop:conformal} 
is canonical in the sense that
it is induced from an extrinsic conformal change of the metric
in the lightcone. If a GCF-manifold has no singular points,
this coincides with the usual conformal change of the
conformally flat metric.
However, if a GCF-manifold admits singular points,
our conformal change may not preserve the
admissibility in general, since the dual metric
$\check g$ may diverge at a degenerate 
point of $g$ (see Remark \ref{rmk:addmissible}).
We also remark that the 
singular sets may not be stable under conformal changes.

\section{Frontal bundles}
\label{sec:frontal}
Let $M^n$ be an  oriented $n$-manifold $(n\ge 1)$ and 
$(M^n,\E,\inner{~}{~},D,\phi)$  a coherent tangent bundle over $M^n$.
Let $\psi:TM^n\to \E$ be another bundle homomorphism
satisfying the following conditions
\begin{enumerate}
 \item $(M^n,\E,\inner{~}{~},D,\psi)$
       is also a coherent tangent bundle,
 \item the pair $(\phi,\psi)$ of bundle 
       homomorphisms satisfies a compatibility condition
       \begin{equation}\label{eq:compati}
	\inner{\phi(X)}{\psi(Y)}=\inner{\phi(Y)}{\psi(X)},
       \end{equation}
       where $X,Y\in T_pM^n$ and $p\in M^n$.
\end{enumerate}
Then $(M^n,\E,\inner{~}{~},D,\phi,\psi)$ is called a {\em frontal
bundle} (see \cite{SUY5}).
The bundle homomorphisms $\phi$ and $\psi$ are called
the {\em first homomorphism\/} and the {\em second homomorphism},
respectively.
We set
\begin{align*}
 \first(X,Y):=&\inner{\phi(X)}{\phi(Y)}, \\
 \second(X,Y):=&-\inner{\phi(X)}{\psi(Y)},\\
 \third(X,Y):= &\inner{\psi(X)}{\psi(Y)} 
\end{align*}
for $X,Y\in T_pM^n$ ($p\in M^n$),  and call them 
{\em the first, the second and the third fundamental forms},
respectively. 
They are all symmetric covariant tensors on $M^n$.
A frontal bundle $(M^n,\E,\inner{~}{~},D, \phi,\psi)$
is called  a {\it front bundle\/} if
\begin{equation}\label{eq:front}
  \Ker(\phi_p)\cap \Ker(\psi_p)=\{0\}
\end{equation}
holds for each   $p\in M^n$.
The conditions for $\phi$ and $\psi$ in the
definition of frontal bundles are symmetric in $\phi$ and
$\psi$, so we can exchange their roles.
(Then the first fundamental form becomes the third fundamental form.)
\begin{ex}\label{ex:front-bundle}
 Let $\bigl(\tilde N^{n+1}(c),\tilde g\bigr)$ be 
 the $(n+1)$-dimensional $1$-connected 
 space form of constant curvature $c$,
 and denote by $\tilde \nabla$ the Levi-Civita connection 
 of $\tilde N^{n+1}(c)$.
 Let $f:M^n\to \tilde N^{n+1}(c)$ be a co-orientable frontal,
 that is, the unit normal vector field $\nu$ is defined globally
 on $M^n$.
 Since the coherent tangent bundle $\E_f$ given in
 Example~\ref{ex:front} is orthogonal to $\nu$,
 we can define a bundle homomorphism
 \[
    \psi_f:T_pM^n\ni X\longmapsto \tilde\nabla_X\nu \in (\E_f)_p
            \qquad (p\in M^n),
 \]
 which can be considered as the shape operator of $f$. 
 Then $(M^n,\E_f,\inner{~}{~},D,\phi_f,\psi_f)$ is a 
 frontal bundle  {\rm(}see  Fact~\ref{fact:frontal} later{\rm)}.
 Moreover, this is a front bundle
 if and only if $f$ is a front,
 which is equivalent to $\first+\third$ being positive definite.
\end{ex}

\begin{fact}[\cite{SUY5}]
\label{fact:frontal}
 Let $f: M^n\to \tilde N^{n+1}(c)$ be a  co-orientable  frontal,
 and $\nu$ a unit normal vector field.
 Then 
 $(M^n,\E_f,\inner{~}{~}, D, \phi_f,\psi_f)$
 as in Example~\ref{ex:front-bundle} is a frontal bundle.
 Moreover, the following identity 
 {\rm (}i.e. the Gauss equation{\rm)}
 holds{\rm:}
 \begin{multline}\label{eq:G}
  \langle R^D(X,Y)v,w\rangle=
  c
  \det\pmt{%
    \inner{\phi(Y)}{v}& \inner{\phi(Y)}{w} \\
    \inner{\phi(X)}{v}& \inner{\phi(X)}{w} 
  } \\
  +
  \det\pmt{%
    \inner{\psi(Y)}{v}& \inner{\psi(Y)}{w} \\
    \inner{\psi(X)}{v}& \inner{\psi(X)}{w} 
  },
 \end{multline}
 where $\phi=\phi_f$ and $\psi=\psi_f$, 
 $X$ and $Y$ are vector fields on $M^n$,
 $v$ and $w$ are sections of $\E_f$,
 and $R^D$ is the curvature tensor of the connection $D${\rm:}
 \[
    R^D(X,Y)v:=D_XD_Yv-D_YD_Xv-D_{[X,Y]}v.
 \]
 Furthermore, this frontal bundle 
 is a front bundle if and only if $f$ is a front.
\end{fact}
Two frontal bundles over $M^n$  
are {\em isomorphic\/} or {\em equivalent\/} if 
there exists an orientation preserving  bundle isomorphism 
between them which preserves the inner products, 
the connections and the bundle maps.
%Let $\bigl(\tilde N^{n+1}(c),g\bigr)$ be the $(n+1)$-dimensional
%$1$-connected space form of constant curvature $c$.
\begin{fact}[\cite{SUY5}]%
\label{thm:fundamental}
 Let  $(M^n,\E,\inner{~}{~},D,\phi,\psi)$ be a 
 frontal bundle over a $1$-connected manifold $M^n$
 $(n\ge 1)$ satisfying \eqref{eq:G},  where $c$ is a real number.
 Then there exists a frontal
 $f\colon{}M^n \to \tilde{N}^{n+1}(c)$
 such that $\E$ is isomorphic to 
 $\E_{f}$  induced from $f$ as in Fact~\ref{fact:frontal}.
 Moreover, such an $f$ is unique up to orientation preserving isometries
 of $\tilde N^{n+1}(c)$.
\end{fact}
Recall that $Q^{n+1}_+$ (resp.\ $Q^{n+1}_-$) is 
the upper (resp.\ lower) lightcone in 
Lorentz-Minkowski space $\bigl(\Lor^{n+2},\inner{~}{~}\bigr)$
of signature $(-+\cdots +)$:
\[
   Q^{n+1}_{\pm} =
   \bigl\{
   z=(z^0,z^1,\dots,z^{n+1})\in\Lor^{n+2}\,;\,
   \inner{z}{z}=0,~\pm z^0>0
   \bigr\}.
\]
From now on, we shall apply Fact \ref{thm:fundamental} to
hypersurface theory in the lightcone $Q^{n+1}_+$: 
A $C^\infty$-map $x:M^n\to Q^{n+1}_+$ 
is called a {\it spacelike frontal\/} if
there exists another $C^\infty$-map $y:M^n\to Q^{n+1}_-$
such that 
\begin{align}\label{eq:ord0}
 & \inner{x}{x}=\inner{y}{y}=0,\quad \inner{x}{y}=1,\\
 & \inner{dx}{y}=\inner{dy}{x}=0,\label{eq:ord1}
\end{align}
where $\inner{dx}{y}$ and $\inner{dy}{x}$ are considered as $1$-forms,
for example, $\inner{dx}{y}$ is defined by
$TM^n\ni X\mapsto \inner{dx(X)}{y}\in \R$.

In this setting, $y$ is called the {\em dual\/} of $x$.
Then $y$ is also a frontal.
Moreover, if the pair 
\[
   (x,y):M^n\longrightarrow Q^{n+1}_+\times Q^{n+1}_-
\]
gives an immersion, $x$ is called a {\it spacelike front}.

\begin{rmk}
 Let $x:M^n\to Q^{n+1}_+$ be a spacelike frontal as above.
 Then the linear map
 \[
   L_p:T_{x(p)}Q^{n+1}_+\ni v \longmapsto \inner{y(p)}{v}\in \R
 \]
 induces a Legendrian lift
 \[
     [L]:M^n\ni p \longmapsto [L_p]\in P(T^*Q^{n+1}_+)
 \]
 of $x$, 
 where $T^*Q^{n+1}_+\ni \alpha\mapsto [\alpha]\in P(T^*Q^{n+1}_+)$
 is the canonical projection.
 Thus, a spacelike frontal is a frontal
 as in Example~\ref{ex:front}.
 Moreover, since
 \[
   \{x(p),y(p)\}^\perp=\{x(p)-y(p),x(p)+y(p)\}^\perp
 \]
 and $x(p)-y(p)$ is a timelike vector,
 the kernel of $L_p$ is a spacelike vector space
 for each $p\in M^n$.

 Conversely, a frontal in $Q^{n+1}_+$
 is spacelike if and only if 
 it has a Legendrian lift
 $[L]:M^n\ni p \mapsto [L_p]\in P(T^*Q^{n+1}_+)$,
 such that the kernel of $L_p$ is a 
 spacelike vector space for each $p\in M^n$.
 In fact, if the kernel $Z_p(\subset T_{x(p)}Q^{n+1}_+)$ 
 of $L_p$ is a spacelike vector space, 
 then the orthogonal complement $(Z_p)^\perp$ in $\Lor^{n+2}$
 is a Lorentzian plane containing $x(p)$.
 Thus, there exists a unique null 
 vector $y(p)\in (Z_p)^\perp$ such that
 $\inner{x(p)}{y(p)}=1$, 
 which induces a $C^\infty$-map $y:M^n\to Q^{n+1}_-$
 and $x$ is a spacelike frontal.
\end{rmk}

\begin{rmk}
 There is a spacelike frontal which is not a front.
 In fact, we set
 \[
     M^n:=S^n\setminus\{(0,\dots,0,\pm 1)\},
 \]
 where
 \[
         S^n=\{(u^1,\dots,u^{n+1})\in \R^{n+1}\,;\, 
        (u^1)^2+\cdots+(u^{n+1})^2=1\},
 \]
 and set
{\small \[
    x:M^n\ni (u^1,\dots,u^{n+1})\longmapsto 
         \frac{1}{\sqrt{2}}\left(
            1,\frac{u^1}{\sqrt{1-(u^{n+1})^2}},\dots,
              \frac{u^{n}}{\sqrt{1-(u^{n+1})^2}},0\right)\in Q^{n+1}_+.
 \]}%
 Then
{\small \[
    y:M^n
       \ni (u^1,\dots,u^{n+1})\longmapsto 
         \frac{1}{\sqrt{2}}
         \left(
            -1,\frac{u^1}{\sqrt{1-(u^{n+1})^2}},\dots,
              \frac{u^{n}}{\sqrt{1-(u^{n+1})^2}},0
         \right)\in Q^{n+1}_-
 \]}%
 gives the dual of $x$. 
 Thus $x$ is a frontal, but not front,
 since % it does not contain the parameter $x_{n+1}$.
 the image of $(x,y)$ lies on an $(n-1)$-dimensional submanifold
 of $Q^{n+1}_+\times Q^{n+1}_-$.
 On the other hand, there is a spacelike front which is not
 an immersion (see Corollary~\ref{cor:added}). 
\end{rmk}

We consider two canonical projections
\[
  \Pi_\pm :Q^{n+1}_\pm \lto S^n_{\pm}:=\{z^0=\pm1\}\cap Q^{n+1}_\pm,
\]
where $(z^0,z^1,\dots,z^{n+1})$ is the canonical
coordinate system of $\Lor^{n+2}$, and $S^n_+$ 
(resp. $S^n_-$) is the
sphere embedded in $Q^{n+1}_+$ (resp. in $Q^{n+1}_-$).
We set
\begin{equation}\label{eq:Gpm}
    G_+:=\Pi_+ \circ x,\qquad G_-:=\Pi_- \circ y,
\end{equation}
which are called the {\em Gauss maps} of $x$ and $y$, respectively.
In this setting, the following assertion
can be proved immediately:
\begin{lem}\label{lem:gauss}
 Let $x:M^n\to Q^{n+1}_+$ $(n\ge 1)$ be a spacelike frontal.
 Then $\inner{dx}{dx}$ is non-degenerate 
 if and only if $G_+$ is an immersion.
 In particular, $x$ itself is an immersion.
\end{lem}

\begin{rmk}
 There is an immersed hypersurface in $Q^{n+1}_+$, 
 which is not spacelike.
 For example, the sub-lightcone
 \[
   \{(z^0,z^1,\dots,z^{n+1})\in Q^{n+1}_+\,; z^{n+1}=0\}
 \]
 is an immersed hypersurface but not spacelike.
 Its Gauss map $G_+$ degenerates everywhere on it.
\end{rmk}

Let $x:M^n\to Q^{n+1}_+$ $(n\ge 1)$ be a spacelike frontal.
Then
\[
  f:=\frac{1}{\sqrt{2}}(x-y):
  M^n\longrightarrow H^{n+1}=\tilde N^{n+1}(-1)\subset\Lor^{n+2}
\]
gives a frontal with the unit normal vector field
\[
    \nu:=\frac{1}{\sqrt{2}}(x+y):
    M^n \longrightarrow S^{n+1}_1\subset \Lor^{n+2},
\] 
where $H^{n+1}$ is the hyperbolic $(n+1)$-space and
$S^{n+1}_1$ is the de Sitter $(n+1)$-space
(i.e. the simply connected
 complete Lorentzian space form of constant sectional
 curvature $1$):
\begin{equation}\label{eq:hyp-des}
\begin{aligned}
   H^{n+1} &=
          \{
             z=(z^0,z^1,\dots,z^{n+1})\in\Lor^{n+2}\,;\,
             \inner{z}{z}=-1,~z^0>0
          \},\\
    S^{n+1}_1 &=
          \{
             z=(z^0,z^1,\dots,z^{n+1})\in\Lor^{n+2}\,;\,
             \inner{z}{z}=1
          \}.
\end{aligned}
\end{equation}
Here, $G_+$ and $G_-$ as in \eqref{eq:Gpm}
are called the {\it hyperbolic Gauss maps} of $f$, 
see \cite{I}, \cite{BIR} and \cite{EGM}.
We set
\[
    \xi:= dx = \frac{1}{\sqrt{2}}(df+d\nu),\qquad
\zeta:= dy = -\frac{1}{\sqrt{2}}(df-d\nu).
\]
Since $x$ is a spacelike frontal, we get a frontal bundle
$(M^n,\E,\inner{~}{~},D,\xi,\zeta)$, where $D$ is 
a metric connection of $\E$ induced by the canonical
connection in $\Lor^{n+2}$ by taking the tangential components.
Moreover,  the following assertion holds:
\begin{thm}%
\label{thm:fundamental2}
 Let $M^n$ $(n\ge 1)$ be a $1$-connected $n$-dimensional
 manifold and  $(M^n,\E,\inner{~}{~},D,\xi,\zeta)$ a
 frontal bundle satisfying 
 \begin{multline}\label{eq:G2}
  \langle R^D(X,Y)v,w\rangle=
  \det\pmt{%
    \inner{\xi(Y)}{v}& \inner{\xi(Y)}{w} \\
    \inner{\zeta(X)}{v}& \inner{\zeta(X)}{w} 
  }
 \\
  +
  \det\pmt{%
    \inner{\zeta(Y)}{v}& \inner{\zeta(Y)}{w} \\
    \inner{\xi(X)}{v}& \inner{\xi(X)}{w} 
  },
 \end{multline}
 where  $X$ and $Y$ are vector fields on $M^n$,
 and $v$, $w$ are sections of $\E$.
 Then there exists a spacelike frontal
 \[
    x\colon{}M^n \longrightarrow Q^{n+1}_+
 \]
 with its dual $y$ such that $\inner{dx}{dx}$, 
 $-\inner{dx}{dy}$ and $\inner{dy}{dy}$ are  
 the first, the second, and the third fundamental forms
 of $(M^n,\E,\inner{~}{~},D,\xi,\zeta)$, 
 respectively.
 Conversely, any spacelike frontal of $M^n$ into $Q^{n+1}_+$ 
 is given  in this manner.
\end{thm}

\begin{proof}
 We set
 \[
    \phi:=\frac{1}{\sqrt{2}}(\xi-\zeta),\qquad
    \psi:=\frac{1}{\sqrt{2}}(\xi+\zeta).
 \]
 Then one can easily check that 
 $(M^n,\E,\inner{~}{~},D,\phi,\psi)$ satisfies the conditions of
 Fact \ref{thm:fundamental} for $c=-1$.
 In fact, \eqref{eq:G} with $c=-1$ is equivalent to \eqref{eq:G2}.
 Since $M^n$ is simply connected,
 Fact \ref{thm:fundamental} and a standard
 continuation argument imply that there is 
 a frontal $f:M^n\to H^{n+1}$ with unit normal
 vector field $\nu:M^n\to S^{n+1}_1$.
 We now set
 \[
    x:=\frac1{\sqrt{2}}(f+\nu),\qquad y:=-\frac1{\sqrt{2}}(f-\nu).
 \]
 Then $x$ (resp.\ $y$) is a  map into $Q^{n+1}_+$ (resp.\ $Q^{n+1}_-$)
 and it can be easily checked that 
 $\inner{dx}{dx}$, $-\inner{dx}{dy}$ and $\inner{dy}{dy}$ are equal to
 the first, the second, and the third fundamental forms
 of $(M^n,\E,\inner{~}{~},D,\xi,\zeta)$.
 Conversely, let $x$ be a spacelike frontal in $Q^{n+1}_+$,
 then there exists a dual frontal $y$ 
(cf.\ \eqref{eq:ord0} and \eqref{eq:ord1})
 such that
 \[
    f:=\frac{1}{\sqrt{2}}(x-y),\qquad
    \nu:=\frac{1}{\sqrt{2}}(x+y).
 \]
 Then $f:M^n\to H^{n+1}$ is a frontal, and
 $\nu$ is its unit normal vector field.
 As seen in Example \ref{ex:front-bundle},
 $f$ induces a frontal bundle 
 $(M^n, \E_f,\inner{~}{~},D,\phi_f,\psi_f)$.
 If we set
 \[
   \xi:=\frac1{\sqrt{2}}(\phi_f-\psi_f),
   \qquad \zeta:=-\frac1{\sqrt{2}}(\phi_f+\psi_f),
 \]
 then $(M^n,\E,\inner{~}{~},D,\xi,\zeta)$ is the desired
 frontal bundle satisfying \eqref{eq:G2}.
\end{proof}

\section{Conformally flat manifolds and hypersurfaces in $H^{n+1}$}
\label{sec:hyperbolic}
Theorem \ref{thm:main1} follows immediately from the following
\begin{lem}\label{lem:9}
 Let  $(M^n,\E,\inner{~}{~},D,\xi,\zeta)$ be a 
 frontal bundle over a $1$-connected manifold $M^n$
 $(n\ge 3)$  satisfying \eqref{eq:G2}.
 Suppose that the regular set $\Reg_{M^n}$ 
 of $\xi$ is dense in $M^n$. 
 Then the first fundamental form $g:=\xi^*\inner{~}{~}$
 is a conformally flat metric on $\Reg_{M^n}$, and the
 Schouten tensor $A$ of $g$ coincides with $-\second$,
 where $\second$ is 
 the second fundamental form of the frontal bundle.
 Moreover, the third fundamental form 
 $\check g:=\zeta^*\inner{~}{~}$ satisfies
 \eqref{eq:tilde}.
\end{lem}

\begin{proof}
 The curvature tensor $R_g$ of $g$ is
 related to $R^D$ by
 \[
    g(R_g(X,Y)Z,W)=\inner{R^D(X,Y)\xi(Z)}{\xi(W)},
 \]
 where $X,Y,Z,W$ are vector fields on $M^n$.
 Substituting \eqref{eq:G2} to $v=\xi(Z)$ and $w=\xi(W)$, and by
 contraction,
 the Ricci tensor $\Ric_g$ is given by
 \begin{equation}\label{eq:ricci}
   \Ric_g=-\trace_{\first}(\second) g-(n-2)\second,
 \end{equation}
 where $\trace_{\first}$ denotes the trace with respect to 
 the first fundamental form $g=I=\xi^*\inner{~}{~}$.
 Then the scalar curvature $S_g$ is given by
 \begin{equation}\label{eq:scal}
     S_g=-2(n-1)\trace_{\first}(\second).
 \end{equation}
 When $n=2$, this implies the equivalency between the
 Gaussian curvature and the mean curvature mentioned in the
 introduction.
 On the other hand, when 
 $n\ge 3$, by \eqref{eq:ricci} and \eqref{eq:scal},
 we have that
 \[
   -\second =\frac{1}{n-2}\left(\Ric_g-\frac{S_g}{2(n-1)}g\right)=A,
 \]
 where $A$ is the Schouten tensor as in \eqref{eq:schouten}.
 Then $A$ is a Codazzi tensor because of \eqref{eq:c} for $\xi$ and
 $\eta$.
 Moreover, if $n\geq 4$,
 one can easily see that the equation \eqref{eq:G2} with $v=\xi(Z)$
 and $w=\xi(W)$ is equivalent to having that
 the conformal curvature tensor as in \eqref{eq:weyl}
 vanishes identically  on $\Reg_{M^n}$.
 Let $(u^1,\dots,u^n)$ be a local coordinate system 
 of $\Reg_{M^n}$.
 We set 
 \[
   dy(\partial_i)=\sum_j \lambda_i^j dx(\partial_j)+c_1 x+ c_2 y,
 \]
 where $\partial_i=\partial/\partial u^i$.
 Since $\inner{x}{dy}=\inner{x}{dx}=0$ and
 $\inner{y}{dx}=\inner{y}{dy}=0$, we have that
 $c_1=c_2=0$ and
 \[
    A_{ij}=\inner{dx(\partial_i)}{dy(\partial_j)}
        =\sum_{k} \lambda_i^k g_{kj}.
 \]
 Then it holds that
 $\lambda_i^j=\sum_{k}A_{ik}g^{kj}$ and
 \begin{align*}
  \check g_{ij}=\inner{dy(\partial_i)}{dy(\partial_j)} 
     =\sum_{a,b} \lambda_i^a \lambda_j^b g_{ab}
      =\sum_{a,b} A_{ia}g^{ab}A_{bj},
 \end{align*}
 which proves the assertion, where $(g^{ab})$ is the
 inverse matrix of $(g_{ij})$.
\end{proof}
\begin{proof}[Proof of Theorem \ref{thm:main2}]
 Let $(M^n,\E,\inner{~}{~},D,\phi)$ be an admissible GCF-manifold.
 By \eqref{eq:tilde-phi}, $\check \phi$ gives
 a bundle homomorphism between $TM^n$ and $\E$.
 By the previous Lemma~\ref{lem:9}, 
 $\check\phi^*\inner{~}{~}$
 is a conformally flat metric on $\Reg^*_{M^n}$.
 By Theorem~\ref{thm:main1},
 the Schouten tensor $A$ is common in two metrics $\phi^*\inner{~}{~}$ 
 and $\check\phi^*\inner{~}{~}$.
 As pointed out in Section~\ref{sec:duality}, 
 the second Bianchi identity with  respect to $\phi^*\inner{~}{~}$
 implies that $A$ is a Codazzi tensor on $\Reg^*_M$, and then
 it is equivalent to the relation
 \begin{equation}\label{eq:tildephi}
    \nabla_X \hat A(Y)-\nabla_Y\hat A(X)-
\hat A([X,Y])=0,
 \end{equation}
 where $\hat A$ is given in \eqref{hatA} and
 $\nabla$ is the Levi-Civita connection of
 $\phi^*\inner{~}{~}$.
 Since 
 $(\Reg^*_{M^n}, TM^n|_{\Reg^*_{M^n}}, \phi^*\inner{~}{~}, \nabla, \hat A)$
 is isomorphic to
 $(\Reg^*_{M^n},\E|_{\Reg^*_{M^n}}, \inner{~}{~}, D,\check \phi)$
 by $\phi$,
 the identity  \eqref{eq:tildephi} yields that
 \begin{equation}\label{eq:tildephi2}
   D_X\check\phi(Y)-D_Y\check \phi(X)-
       \check\phi([X,Y])=0
 \end{equation}
 holds on $\Reg^*_{M^n}$.
 Since $\Reg^*_{M^n}$ is dense in $M^n$,
 \eqref{eq:tildephi2} holds on the whole of $M^n$.
 Thus 
 $(M^n,\E, \inner{~}{~},D,\check \phi)$
 is also a GCF-manifold.
 Since the set of regular points $\Reg_{M^n,\check\phi}$ of
 $\check\phi$ contains $\Reg^*_{M^n}$, 
 $(M^n,\E, \inner{~}{~},D,\check\phi)$ is an admissible
 GCF-manifold.
 Then as seen in the proof of Lemma \ref{lem:9},
 $(M^n,\E,\inner{~}{~},D,\phi,\check\phi)$ is a 
 frontal bundle over $M^n$  satisfying \eqref{eq:G2} for
 $\xi=\phi$ and $\zeta=\check\phi$.
 Since the condition \eqref{eq:G2}
 is symmetric with respect to $\xi$ and $\zeta$,
 by switching the roles of $\phi$ and $\check \phi$,
 we can conclude that
 $\check{\check \phi}$ coincides with $\phi$.
\end{proof}

Let $(M^n,\E,\inner{~}{~},D,\phi)$ be an admissible GCF-manifold. 
Then as seen in the above proof of Theorem \ref{thm:main2}, 
it induces a bundle homomorphism $\check\phi:TM^n\to \E$
satisfying \eqref{eq:G2} for
$\xi=\phi$ and $\zeta=\check \phi$.
As seen in the proof of Theorem \ref{thm:fundamental2},
if we set
\[
  \psi_1:=\frac{1}{\sqrt{2}}(\phi-\check \phi),\qquad
  \psi_2:=\frac{1}{\sqrt{2}}(\phi+\check \phi),
\]
then $(M^n,\E,\inner{~}{~},D,\psi_1,\psi_2)$ 
satisfies the condition of Fact \ref{thm:fundamental}
for $c=-1$. Then we get a frontal
$f:\tilde M^n \to H^{n+1}$,
where $\tilde{M} ^n$ is the universal covering of $M^n$.
Here, $f$ is determined up to rigid motions in $H^{n+1}$.
Let $\nu$ be the unit normal vector field of $f$.
Then it induces a spacelike wave front (cf. \cite[Section 2]{SUY5})
$\nu:M^n\to S^{n+1}_1$,
where $S^{n+1}_1$ is the de Sitter space.
The two frontals $f$ and $\nu$ are 
mutually dual objects. The realization theorem of
spacelike fronts in $S^{n+1}_1$ is proved in \cite{SUY5} using the
duality. In particular, if $M^n$ is $1$-connected,
we get two maps
\begin{align*}
   \Psi&:\GCF(M^n)\lto \MFr(M^n,H^{n+1}),\\
   \Psi^*&:\GCF(M^n)\lto \MFr(M^n,S^{n+1}_1),
\end{align*}
at the same time,
where $\GCF(M^n)$ is the set of admissible GCF-manifolds modulo
structure-preserving bundle isomorphisms 
and the set $\MFr(M^n,H^{n+1})$ (resp. $\MFr(M^n,S^{n+1}_1)$) 
is the set of congruent classes of frontals in $H^{n+1}$ 
(resp.\ spacelike frontals in $S^{n+1}_1$).
The maps $\Psi$ and $\Psi^*$
are both two-fold maps, since $\Psi$ 
(resp.\ $\Psi^*$) 
takes the same values for a dual GCF-manifold as 
for the given GCF-manifold.

We fix an admissible GCF-manifold
$(M^n,\E,\inner{~}{~},D,\phi)$. 
By Lemma \ref{lem:gauss}, the intersection of
two regular sets of the hyperbolic Gauss maps 
$G_\pm:M^n \to S^{n}$ of $f$ coincides with
the set $\Reg_{M^n}^*$. 
Let $\MFr^*(M^n,H^{n+1})$
(resp. $\MFr^*(M^n,S^{n+1}_1)$) 
the subset of $\MFr(M^n,H^{n+1})$
(resp. $\MFr(M^n,S^{n+1}_1)$)
consists of wave fronts 
(resp.\ spacelike wave fronts)
whose pairs of  Gauss maps $(G_+,G_-)$
both have dense regular sets in $M^n$.
Thus we get the following
\begin{thm}
 Let $M^n$ $(n\ge 3)$ be a $1$-connected manifold.
 Suppose that $\GCF(M^n)$ is non-empty.
 Then 
 $\Psi:\GCF(M^n)\to \MFr^*(M^n,H^{n+1})$
 {\rm(}resp.\  $\Psi^*:\GCF(M^n)\to 
 \MFr^*(M^n,S^{n+1}_1)${\rm)}
 is a surjective two-fold map.
 An admissible GCF-manifold has no singular points
 if and only if the positive Gauss map
 $G_+:M^n\to S^n$ is an immersion.
\end{thm}

\begin{rmk}
 An embedding $f:S^n\to H^{n+1}$ is called 
 {\em horo-regular\/} if at least one of the hyperbolic Gauss maps of $f$
 is a diffeomorphism.
 On the other hand, if $f$ lies in the closure of
 the interior of the osculating horosphere, $f$ is called 
 {\it horoconvex\/} (cf. \cite{BIR}).
 A horo-regular horoconvex embedding
 $f:S^n\to H^{n+1}$ is called {\it strictly horoconvex}.
 In \cite{EGM}, it is pointed out that 
 an embedding $f:S^n\to H^{n+1}$ is strictly horoconvex if and only if
 both Gauss maps $G_+$ and $G_-$ are diffeomorphisms.
 Several characterizations of horoconvexity are given in
 \cite{BIR} and \cite{EGM}.
 When $f\in \MFr^*(M^n,H^{n+1})$ 
 is a horo-regular embedding, \cite{EGM} showed that
 there is a conformally flat metric $g$ on $S^n$ realizing $f$. 
 Our map $\Psi$ is a generalization of this procedure in \cite {EGM}.
 In \cite{EGM}, horo-regularity (resp.\ strict horoconvexity)
 is called horospherical convexity (resp.\ strongly $H$-convexity).
 When $f$ is a front in $H^{n+1}$, then parallel family of wave front
 $\{f_\delta\}_{\delta\in \R}$ is induced.
 Like as in the case of horo-regular hypersurfaces
 in \cite {EGM}, $f_\delta$ induces an admissible GCF-manifold
 whose metric is a scalar multiple of the metric
 $\inner{~}{~}$ of the GCF-manifold induced by $f$.
\end{rmk}

\begin{cor}[Kuiper \cite{K}]\label{cor:Kuiper}
 Let $M^n$ $(n\ge 3)$ be a compact $1$-connected manifold.
 Then $M^n$ admits a conformally flat metric if and only if
 $M^n$ is diffeomorphic to $S^n$.
\end{cor}

\begin{proof}
 Suppose that there is a conformally flat metric $g$
 on $M^n$, then $G_+:M^n\to S^n$ is an immersion.
 Since $M^n$ is compact, it gives a finite covering map.
 Since $M^n$ and $S^n$ are both $1$-connected,
 $G_+$ must be bijective.
\end{proof}

Contrary to the above assertion, we can prove the following
\begin{prop}\label{prop:admissible}
 There is a compact $1$-connected 
 admissible generalized  conformally 
 flat manifold that is not homeomorphic to a 
 sphere.
\end{prop}
\begin{proof}
 Take a generalized Clifford torus $S^2\times S^{n-2}$ in $S^{n+1}$.
 By a conformal transformation, we can get a
 hypersurface immersed in an open hemisphere which can be
 identified with the hyperbolic space $H^{n+1}$.
 Then, we get an immersion $f:S^2\times S^{n-2}\to H^{n+1}$
 with a unit normal vector field $\nu$. 
 Then $x:=(f-\nu)/\sqrt{2}$
 gives a front in $Q^{n+1}_+$ and 
 the metric $\inner{dx}{dx}$ induces the desired 
 generalized conformally flat structure
 on $S^2\times S^{n-2}$.
\end{proof}

\begin{cor}\label{cor:added}
 There is a spacelike front in $Q^{n+1}_+$ which is not
 an immersion.
\end{cor}

\begin{proof}
 Let $x:S^{2}\times S^{n-2}\to Q^{n+1}_+$ be
 the spacelike front as in the proof of 
 Proposition~\ref{prop:admissible}.
 If $x$ is an immersion, then the corresponding 
 compact $1$-connected 
 admissible generalized  conformally 
 flat manifold has no singularity.
 Then by Kuiper's theorem (i.e. Corollary \ref{cor:Kuiper}),
 it is diffeomorphic to $S^n$, which makes a contradiction.
\end{proof}

At the end of this section, we discuss conformal changes 
of a given front in $Q^{n+1}_+$.
\begin{proof}[Proof of Proposition \ref{prop:conformal}]
 Since the assertion is a local property, we may assume that
 $M^n$ is simply connected.
 Then by Lemma~\ref{lem:9}, there exists a spacelike frontal
 $x:M^n\to Q^{n+1}_+$
 which induces $(M^n,\E,\inner{}{},D,\phi)$.
 We denote by
 $y:M^n\to Q^{n+1}_-$
 its dual. 
 Then 
 \[
   \tilde x:=e^{\sigma}x:U\lto Q^{n+1}_+
 \]
 gives a new immersion whose first fundamental form is given by
 \[
      \tilde g=\inner{d\tilde x}{d\tilde x}
      =e^{2\sigma}\inner{dx}{dx}=e^{2\sigma} g.
 \]
 Let $\tilde\E$ be the subbundle of $\tilde x^*TQ^{n+1}_+$
 perpendicular to $y$, and $\tilde D$ be an induced connection on 
 $\tilde \E$ with respect to the canonical connection 
 on $\Lor^{n+2}$.  Then it  induces a new GCF-manifold.
\end{proof}

Let $U$ be a domain in $S^n$, and
$x:U\to Q^{n+1}_+$ a canonical embedding, that is,
\[
   x=\pmt{1 \\ p}\subset \Lor^{n+2}
\]
and $p:U\to S^n\subset \R^{n+1}$ is the canonical inclusion.
Then 
\[
   \tilde x:=e^{\sigma}x:U\lto Q^{n+1}_+
\]
gives a new immersion whose first fundamental form is given by
\[
      \tilde g=\inner{d\tilde x}{d\tilde x}
      =e^{2\sigma}\inner{dx}{dx}=e^{2\sigma} g,
\]
where $g$ is the induced metric of $U$ from the unit sphere $S^n$.
Recall that the dual of $\tilde x$ is a map
$\tilde y:U\longrightarrow  Q^{n+1}_-$
such that 
\[
  \inner{\tilde x}{\tilde y}=1,\quad 
  \inner{d\tilde x}{\tilde y}
   =\inner{d\tilde y}{\tilde x}=0.
\]
Then one can directly verify that
\begin{align}
 \tilde y&=
   -\frac12\Delta \tilde x -
  \frac{\inner{\Delta \tilde x}{\Delta \tilde x}}8\tilde x
 \label{eq:1}
 \\
  &=\frac{e^{-\sigma}}{2}\left\{ 
    \pmt{-1\\ p}-|d\sigma|^2 \pmt{1 \\ p}-2 \pmt{0\\ \alpha}
 \right\} 
 \label{eq:2}
\end{align}
holds, where
\[
   \Delta:=\sum_{i,j}g^{ij}\nabla_j \nabla_i,\quad 
   |d\sigma|^2:=\sum_{i,j}g^{ij}\sigma_i\sigma_j,\quad
   \alpha:=\sum_{j,k}g^{jk}\sigma_j p_k:U\to \R^{n+1}.
\]
Here, $g_{ij}$ ($i,j=1,\dots,n$) are the components of the 
metric $g$ with respect to a local coordinate system
$(u^1,\dots,u^n)$,
$(g^{ij})$ is the inverse matrix of $(g_{ij})$, 
$\sigma_j=\partial\sigma/\partial u^j$ and $p_j=\partial p/\partial u^j$.
The first equation \eqref{eq:1} is the formula in \cite{LJ}.

In particular, the second fundamental form of $\tilde x$
is given by
\[
  \second=-\inner{d\tilde x}{d\tilde y}=
  \left(\frac{1+|d\sigma|^2}{2}\right)g+d\sigma\otimes d\sigma
   -\Hess(\sigma).
\]
The symmetric covariant tensor $\second$ satisfies the Codazzi equation,
since symmetricity of  $\nabla \second$ 
is equivalent to the condition \eqref{eq:c} for the second homomorphism.

\begin{rmk}
\label{rmk:addmissible}
 The frontal bundle 
 $(M^n,\tilde\E,e^{2\sigma}\inner{~}{~},D,d\tilde x)$
 constructed in the proof of Proposition \ref{prop:conformal}
 might not be admissible, since
 $\tilde y$ given in \eqref{eq:2} can diverge
 if $\tilde x$ admits singular points.
 This generalization of the conformal change of
 Riemannian metrics
 is somewhat related to 
 the hyperbolic Christoffel problem  posed in \cite{EGM}.
\end{rmk}

\section{Isometric immersions of a Riemannian $2$-manifold  into $Q^3_+$}
If $n\ge 3$, any simply connected conformally flat 
Riemannian $n$-manifold is uniquely immersed in the lightcone
$Q^{n+1}_+$.
However, if $n=2$, the situation is different.
In this section, we show that
there is infinite dimensional freedom for isometric 
immersions, 
unless the given simply connected Riemannian $2$-manifold is 
homeomorphic to $S^2$.
More precisely, we can determine the moduli of the set of immersions of a
given simply connected  Riemannian manifold into the $3$-dimensional
lightcone $Q^3_+$.
First, we prove the following:
\begin{prop}\label{thm:embedding}
 Let $(M^2,g)$ be a $1$-connected Riemannian $2$-manifold.
 Then there is an isometric embedding $x:(M^2,g)\to Q^3_+$.
\end{prop}
\begin{proof}
 The well-known uniformization theorem implies that
 $(M^2,g)$ is conformally equivalent to the sphere $S^2$,
 the plane $\C$ or the unit disc $\D^2$.
 Since $\C$ and $\D^2$ are conformally embedded in  
 the unit sphere $S^2$, there is a conformal
 embedding
 \[
    x_1:(M^2,g)\lto S^2.
 \]
 Since $x_1$ is conformal, there is a smooth function 
 $\sigma\in C^\infty(M^2)$
 such that $g=e^{2\sigma}g_{1}$, where $g_1$ is a
 metric of constant Gaussian curvature $1$ 
 induced by $x_1$.
 Then we set
 \[
     x:=e^\sigma x_1:(M^2,g)\lto Q^3_+,
 \]
 which gives the desired isometric embedding.
\end{proof}

On the other hand, as a corollary of Theorem \ref{thm:fundamental2},
the fundamental theorem of surface theory in the lightcone $Q^3_+$
is stated as follows:
\begin{prop}%
\label{thm:fundamental3}
 Let $(M^2,g)$ be a $1$-connected Riemannian manifold
 and  $\second$  a symmetric covariant tensor on $M^2$
 satisfying the Codazzi equation $($i.e.\ the covariant derivative
 $\nabla \second$ with respect to the Levi-Civita
 connection is a symmetric tensor$)$.
 If the Gaussian curvature $K$ of $g$ coincides with 
 the trace of $-\second$, then there exists an isometric
 immersion $x:(M^2,g)\to Q^3_+$ such that
 the second fundamental form of $x$ coincides with $\second$.
 Conversely, any isometric immersions of $(M^2,g)$
 to $Q^3_+$ are given in this manner.
\end{prop}

\begin{proof}
 Let $\nabla$ be the Levi-Civita connection of $g$,
 and $\xi:TM^2\to TM^2$  the identity map.
 We also define a map $\zeta:TM^2\to TM^2$
 so that
 \[
   g\bigl(\xi(X),\zeta(Y)\bigr)=-\second(X,Y).
 \]
 Then $(M^2,T M^2,\inner{~}{~},\nabla,\xi,\zeta)$ is 
 a front bundle.
 In fact, $\zeta$ satisfies \eqref{eq:c} since $\second$ is
 a Codazzi tensor.
 It is sufficient to show the integrability condition
 \eqref{eq:G2} is equivalent to $-K=\trace_{\first}(\second)$.
 To show this,
 we take a local orthonormal frame field
 $e_1$, $e_2$ on $M^2$  such that
 \[
    \lambda_j\delta_{ij}
      =-\second(e_j,e_j)=-\inner{\xi(e_i)}{\zeta(e_j)}
       =-\inner{\xi(e_j)}{\zeta(e_i)}
         \qquad (j=1,2).
 \]
 Then substituting $X=e_1$, $Y=e_2$, $v=e_2$, $w=e_1$,
 \eqref{eq:G2} reduces to
 \[
    K=-\lambda_1-\lambda_2=-\trace_{\first}\second,
 \]
 which proves the assertion.
\end{proof}

\begin{thm}
 Let $(M^2,g)$ be a $1$-connected Riemannian
 manifold,
 and $I_{Q^3_+}(M^2,g)$ the set of congruent classes of
 isometric immersions of $(M^2,g)$ into $Q^3_+$.
 We fix a complex structure of $M^2$ which is
 compatible with respect to the conformal structure.
 Then $I_{Q^3_+}(M^2,g)$ corresponds bijectively
 to the set of holomorphic $2$-differentials on $M^2$
 {\rm (}cf.\ \cite{FK}{\rm)}. 
 In particular,
\begin{enumerate}
 \item if $M^2$ is conformally equivalent to the sphere
       $S^2$, then $I_{Q^3_+}(M^2,g)$ consists of a point,
       namely, the canonical isometric embedding of $(M^2,g)$
       into $Q^3_+$ is rigid,
 \item if $M^2$ is conformally equivalent to the 
       complex plane $\C$, then $I_{Q^3_+}(M^2,g)$ 
       corresponds bijectively to the set $\O(\C)$
       of entire holomorphic functions,
 \item if $M^2$ is conformally equivalent to the 
       unit disc $\D^2$, then $I_{Q^3_+}(M^2,g)$ 
       corresponds bijectively to the set $\O(\D^2)$
       of holomorphic functions on the disc $\D^2$.
\end{enumerate}
\end{thm}
\begin{proof}
 As shown in Proposition \ref{thm:embedding},
 there is an isometric embedding $f:M^2\to Q^3_+$.
 In particular, there is a Codazzi tensor $-\second_K$ 
 whose trace is equal to the Gaussian curvature $K$ of the
 metric $g$.
 Then, by Proposition \ref{thm:fundamental3},
 $I_{Q^3_+}(M^2,g)$ can be identified with
 the set $\Cod_{-K}(M^2,g)$ of Codazzi tensors on 
 $(M^2,g)$
 whose traces are equal to $-K$.
 We denote by 
$\Cod_{0}(M^2,g)$ the set of traceless
 Codazzi tensors on $(M^2,g)$.
 The following map is a bijection;
 \begin{equation}\label{map1}
  \Cod_{0}(M^2,g)\ni B \longmapsto 
     B+\second_K\in \Cod_{-K}(M^2,g).
 \end{equation}
 Since $M^2$ is simply connected, it can be considered
 as a Riemann surface biholomorphic to $S^2$, $\C$ or $\D^2$.
 As shown in \cite{LSW},
 the set of traceless  Codazzi tensors on a compact Riemann surface 
 is identified with the set of holomorphic $2$-differentials on it.
 Thus, if $M^2=S^2$, the vector space  $\Cod_{0}(M^2,g)$
 is zero-dimensional, that is, $f$ is rigid.
 So we consider the remaining cases $M^2=\C$ and $M^2=\D^2$
 with the canonical complex coordinate $z$.
 We denote by $\O(M^2)$ the set of  holomorphic 
functions on $M^2$.
 Then
 \begin{equation}\label{map2}
     \O(M^2)\ni \phi(z) \longmapsto
    \phi(z)dz^2+\overline{\phi(z)dz^2}\in \Cod_{0}(M^2,g)
 \end{equation}
is a linear isomorphism.
 (The set of holomorphic $2$-differentials on $M^2$ can be 
 identified with
 $\O(M^2)$, since there is a globally defined
 holomorphic 2-differential $dz^2$ on $M^2$.)
 Thus, combining the two maps \eqref{map1} and \eqref{map2},
 $\O(M^2)$ can be identified with 
 $I_{Q^3_+}(M^2,g)$.
\end{proof}

Let $M^2$ be a  $2$-manifold and $x:M^2\to Q^3_+$ 
an immersion. 
Then there exists the dual as a $C^\infty$-map 
$y:M^2\to Q^3_-$
such that $\inner{x}{y}=1$ and 
$\inner{dx}{y}=\inner{x}{dy}=0$.
Even if $x$ is an immersion, the dual $y$ of $x$ may have singular
points. 
In fact, $y$ is a spacelike front in $Q^3_-$ in general.
The first author \cite{L} and Izumiya-Saji \cite{IS} pointed out that
if $x$ has zero Gaussian curvature with respect to 
the induced metric,
so does $y$ on its regular set.
It should be remarked that this duality on flat surfaces
corresponds to the following intrinsic duality:

\begin{prop}
 Let $(M^2,g)$ be a flat Riemannian $2$-manifold and
 $\second$ a traceless Codazzi tensor of $(M^2,g)$.
 Then the metric $\check g$ defined by \eqref{eq:tilde}
 is also flat on the regular set of $\check g$.
 Moreover, $\second$ is also a 
 traceless Codazzi tensor with respect to $\check g$ on
 the regular set of $\check g$. 
\end{prop}

\begin{proof}
 This corresponds to the case of $n=2$ of our intrinsic duality.
 (As seen above, the set of traceless Codazzi tensors
 on $(M^2,g)$ can be identified with $I_{Q^3_+}(M^2,g)$
 if $M^2$ is simply connected.)
 The proof is parallel to that of Lemma \ref{lem:9}. 
\end{proof}

\begin{rmk}
 Let $\nabla$ and $\check \nabla$ be the Levi-Civita connections
 of $g$ and $\check  g$ respectively.
 Then $(\second,\check \nabla)$ is the dual Codazzi structure
 of the Codazzi structure $(\second,\nabla)$ in the sense 
of Shima \cite{S}.  
\end{rmk}
Several examples of flat surfaces in $Q^3_+$
are given by the first author \cite{L2}.
Recently, Izumiya and Saji \cite{IS} showed that linear Weingarten 
fronts 
in $Q^3_+$ correspond to linear Weingarten 
fronts in  $H^3$ as a Legendrian duality.
In this setting, one can easily check that
flat fronts in $Q^3_+$
corresponds to flat fronts in $H^3_+$.
On the other hand,
G\'alvez, Mart\'inez and Mil\'an \cite{GMM} found a holomorphic
representation formula for flat surfaces in $H^3$
using a Hessian structure induced by the above
Codazzi structure $(\second,\nabla)$. 
Like as in the case of $n\ge 3$, this duality corresponds
to the following two-fold map
between two sets of flat fronts in $Q^3_+$ and in $H^3$ on $M^2$
\[
   \MFr^0(M^2,Q^{3}_+)\ni (x,y)\longmapsto
     \frac{1}{\sqrt{2}} (x-y,x+y) \in \MFr^0(M^2,H^{3}).
\]
Several examples and global properties of flat fronts 
in $H^3$ are given 
in \cite{KUY}, \cite{KRUY} and \cite{KRUY2}.
\begin{acknowledgement}
The second and the third authors thank Pablo Mira
for fruitful discussions on the occasion of 
his visit at Osaka University in 2009, 
and especially for informing them of  the work \cite{EGM}.
The second author thanks Shyuichi Izumiya
for valuable comments on the first draft of this paper. 
\end{acknowledgement}


\begin{thebibliography}{99}

\bibitem{AD} 
  A.~R.~C.~Asperti and M.~Dajczer,
  {\itshape 
	Conformally flat Riemannian manifolds
	as hypersurfaces of the light cone},
	Canad.\ Math.\ Bull. {\bfseries 32} (1989), 281--285.

\bibitem{Br}
  H.~W.~Brinkmann, 
  {\itshape 
	On Riemann spaces conformal to Euclidean space},
	Proc.\ Nat.\ Acad.\ Sci {\bfseries 9} (1923), 1--3.

\bibitem{BIR}
  M.~Buosi, S.~Izumiya and M.~A.~Ruas,
  {\itshape 
	Horo-tight spheres in hyperbolic space},
	preprint. 

\bibitem{EGM}
  J. M. Espinar, J. A. G\'alvez and P. Mira,
  {\itshape 
	Hypersurfaces in $H^{n+1}$ and
	conformally invariant equations:
	The generalized Christoffel and Nirenberg problems},
	J. Eur. Math. Soc. {\bf 11} (2009), 903-939.

\bibitem{FK}
 H. M. Farkas and I. Kra, 
  {\sc Riemann surfaces} (2nd edition), 
	Springer-Verlag (1992).

\bibitem{GMM}
  J. A. G\'alvez, A. Mart\'inez and F. Mil\'an, 
  {\itshape 
	Flat surfaces in hyperbolic $3$-space}, 
	Math.\ Ann. {\bf 316} (2000), 419-435.

\bibitem{I}
  S.~Izumiya, 
  {\itshape 
	Legendrian dualities and 
	spacelike hypersurfaces in the lightcone},
	Moscow Mathematical Journal {\bf 9} (2009), 325--357.

\bibitem{IS}
 S.~Izumiya and K.~Saji, 
 {\itshape 
	The mandala of Legendrian
	dualities for pseudo-spheres of Lorentz-
	Minkowski space and ``flat'' spacelike surfaces},
	preprint.

\bibitem{K} 
  N.~H.~Kuiper,
  {\itshape 
	On conformally-flat space forms in the large},
	Ann.\ of Math. {\bf 91} (1949), 916-924.

\bibitem{KUY}
  M.~ Kokubu, M.~Umehara, and K.~Yamada, 
  {\itshape
	Flat fronts in hyperbolic $3$-space},
	Pacific J. Math. {\bf 216} (2004), 149--175. 

\bibitem{KRUY}
  M.~ Kokubu, W.~Rossman, M.~Umehara and K.~ Yamada,  
  {\itshape
	Flat fronts in hyperbolic $3$-space and their caustics}, 
	J. Math.\ Soc.\ Japan {\bfseries 59} (2007), 265--299. 

\bibitem{KRUY2}
  M.~Kokubu, W.~Rossman, M.~Umehara and K.~Yamada,  
  {\itshape
	  Asymptotic behavior of flat surfaces in hyperbolic 3-space}, 
	J. Math.\ Soc.\ Japan {\bf 61} (2009), 799-852.

\bibitem{L} 
  H.~L.~Liu,
  {\itshape 
	Surfaces in the light cone},
	J. Math.\ Annal. {\bf 325} (2007), 1171-1181.

\bibitem{L2} 
  H.~L.~Liu,
  {\itshape 
	Maximal surfaces in $3$-dimensional lightlike cone $Q^3$}, 
	preprint.

\bibitem{LJ} 
  H.~L.~Liu and S.~D.~Jung,
  {\itshape 
	Hypersurfaces in lightlike cone},
	J. of Geom.\ and Phys. {\bf 58} (2008), 913--922.

\bibitem{LSW} 
  H.~L.~Liu, U.~Simon and C.~P.~Wang,
  {\itshape 
	Codazzi Tensors and the Topology of Surfaces},
	Ann.\ Global Anal.\ Geom.
	{\bfseries 16} (1998), 189--202.

\bibitem{S} 
  H. Shima, 
  {\sc 
	The Geometry of Hessian Structures}, 
	World Scientific Publishing Company, 2007.  

\bibitem{SUY1} 
  K. Saji, M. Umehara and K. Yamada,
  {\itshape 
	The geometry of fronts}, 
	Ann.\ of Math. {\bf 169} (2009), 491--529.
 
\bibitem{SUY2} 
  K. Saji, M. Umehara and K. Yamada,
   {\itshape 
	Behavior of corank one singular points on wave fronts},
	Kyushu J. of Mathematics {\bf 62} (2008), 259--280.

\bibitem{SUY4} 
  K. Saji, M. Umehara and K. Yamada,
  {\itshape
	The duality between singular points and
	inflection points on wave fronts}, 
	to appear in Osaka J. Math., arXiv:0902.0649. 

\bibitem{SUY5} 
  K. Saji, M. Umehara and K. Yamada,
  {\itshape
	The intrinsic duality of wave fronts},
	preprint, arXiv:0910.3456.

\end{thebibliography}
\end{document}